\documentclass[11pt]{amsart}
\usepackage{graphicx}
\pdfoutput=1

\issueinfo{00}{0}{Month}{Year}
\copyrightinfo{2008}{University of Calgary}
\pagespan{1}{2}
\PII{ISSN 1715-0868}

\newtheorem{thm}{Theorem}[section]
\newtheorem{proposition}[thm]{Proposition}
\newtheorem{lemma}[thm]{Lemma}
\newtheorem{corollary}[thm]{Corollary}
\newtheorem{conjecture}[thm]{Conjecture}

\theoremstyle{definition}

\newcounter{algorithm}
\renewcommand{\thealgorithm}{\thesection.\arabic{algorithm}}
\newenvironment{algorithm}[1]{%
	\null
	\refstepcounter{algorithm}%
	\hrule%
	\vspace{0.2em}%
	\noindent\textbf{Algorithm \thealgorithm} #1
	\vspace{0.2em}%
	\hrule%
	\vspace{0.2em}
}{
\vspace{0.2em}%
\hrule%
\null
}

\begin{document}

\title{Maximum spanning trees in normed planes}

\author{Javier Alonso}
\address{Instituto de Matem\'{a}ticas (IMUEx), UEx, 06006 Badajoz, Spain}
\email{jalonso@unex.es}

\author{Pedro Mart\'{\i}n}
\address{Instituto de Matem\'{a}ticas (IMUEx)-Departamento de Matem\'{a}ticas, UEx, 06006 Badajoz, Spain}
\email{pjimenez@unex.es}

%


\thanks{The second author is partially supported by Junta de Extremadura/FEDER Grants Number IB18023 and  GR18023.}
\begin{abstract}
Extending some properties from the Euclidean plane to any normed plane, we show the validity of the Monma-Paterson-Suri-Yao algorithm for finding the maximum-weighted spanning tree of a set of $n$ points, where the weight of an edge is the distance between the end points measured by the norm and there are not repeated distances. For strictly convex normed planes, we expose an strategy for moving slightly the points of the set in order to obtain distinct distances.
\end{abstract}

\maketitle
	
\section{Introduction and notation}

Let $\mathbb{E}^2$ be the Euclidean plane, and $\mathbb{M}^2$ be
a Minkowski plane, namely, $\mathbb{R}^2$ endowed with a norm $\|\cdot\|$. 
For a set $S$ of $n$ points in $\mathbb{M}^2$,
let us consider the undirected graph $G$, where the vertices are the $n$ points and the edges connect every pair of points and are weighted by the distance between them measured with the norm. A \textit{path} in a graph is a sequence of edges such that all the vertices (and therefore also all the edges) are distinct.
A \emph{tree} 
is a graph in which any two vertices are connected by exactly one path. A \emph{maximum spanning tree} of $G$ is a tree having the maximum total edge weight. An algorithm for computing an Euclidean maximum spanning tree is presented in \cite{MPSY}. The algorithm runs in $O(n\log h)$ time and requires $O(n)$ space, where $h$ denotes the number of vertices of the convex hull of $S$. Maximum spanning trees for $L_1$ and $L_{\infty}$ metrics can be computed in $O(n)$ time and $O(n)$  space (\cite{Gabow}).

The notion of Voronoi diagram is widely used in many areas of science. Klein \cite{Klein} introduces \textit{abstract Voronoi diagram} as a unifying approach of Voronoi diagram centred on bisecting curves rather on distances. The algorithm in \cite{MPSY} is based on the construction of the furthest Voronoi diagram of $S$, that is a subdivision of the plane into regions such that the region corresponding to a point $p\in S$ is the locus of points further away from $p$ than from other point of $S$. The furthest Voronoi diagram can be computed using different algorithms in $O(n\log n)$ time for the Euclidean and  $L_p$ ($1\leq p \leq \infty$) norms (\cite{Au}, \cite{Lee}, \cite{Shamos-Hoey}). Mehlhorn et al. (\cite{Mehlhorn}) provide a randomized algorithm to construct the furthest abstract Voronoi diagram in expected $O(n\log n)$ time.

Following the skeleton of the paper of Monma et al. (\cite{MPSY}), we present in Section \ref{structure} the structure of the unique maximum spanning tree in $\mathbb{M}^2$ for $n$ points that do not present equalities of distances, and  prove in Section \ref{algorihtm} that the approach in (\cite{MPSY}) for computing the maximum spanning tree of a set of $n$ points works also in $\mathbb{M}^2$. The computation takes $O(n\log n)$ time whether a subroutine for the construction of any furthest point Voronoi diagram taking $O(n\log n )$ is available. This is the case of a wide range of norms, for instance $L_p$ ($1<p<\infty$), as we comment in the previous paragraph. In Section \ref{applications} we expose some applications to other geometrical computation problems.
We prove in last section that given $n$ points in a strictly convex plane $\mathbb{M}^2$, and $\epsilon>0$, each point can be moved a distance less than $\epsilon$ in such a way that all the distances between the $n$ points are different. A similar strategy does not work whether the plane is not strictly convex. 

\section{Notation and preliminaries}\label{notation}
\begin{figure}[h]
	\begin{center}
        \scalebox{0.8}{\includegraphics{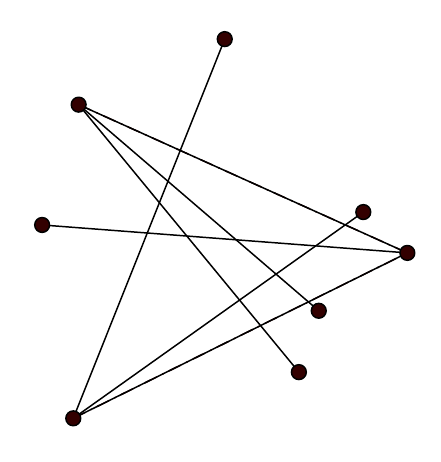}}
\\\caption{\label{mxstree}$\mathtt{MXST}(S)$ of a set $S$ of $n$ points in $\mathbb{E}^2$.}
	\end{center}
\end{figure}

Let $S=\{p_1,p_2,...,p_n\}$ be a set of $n$ points in a normed plane $\mathbb{M}^2$. From Section \ref{notation} to Section \ref{applications} we consider that all the distances between any two points are distinct. This condition guarantees that there exists a unique maximum spanning tree of $S$, denoted by $\mathtt{MXST}(S)$ (see Figure \ref{mxstree}). Given $p_i\in S$, we call $p_k\in S$ the \textit{furthest neighbour} of $p_i$ if $\|p_i-p_k\|=\mathrm{max}\{\|p_i-p_j\|:\ j=1,...,n\}$. We denote $\mathtt{FNG}(S)$ the unique directed graph whose vertices are the points of $S$ and whose directed edges connect each point of $S$ to its furthest neighbour. 

We call $B(x,r)$ the \textit{ball with center $x\in \mathbb{M}^2$ and radius $r>0$}. 
The convex hull of a set $A$ is denoted by $\mathrm{conv}(A)$, and its boundary by $\partial \mathrm{conv}(A)$. Given $x,y\in \mathbb{M}^2$, $\overline{xy}$ denotes the \textit{line segment} connecting the two points, $L(x,y)$ denotes the line through $x$ and $y$, 
and $\mathrm{bis}(x,y)$ denotes the \textit{bisector} for $x$ and $y$, namely, the set of points equidistant from $x$ and $y$. We say that  $x$ is Birkhoff orthogonal to $y$ (denoted by $x\perp y$) if $\|x\|\leq \|x+\lambda y\|$ for every $\lambda\in \mathbb{R}$. Birkhoff orthogonality extends the notion of Euclidean orthogonality to $\mathbb{M}^2$.

\begin{lemma}\label{furconvex}
Let $S$ be a set of $n$ points in $\mathbb{M}^2$. Then,
\begin{enumerate}
\item Given $p\in \mathbb{M}^2$, all the points where the furthest distance from $S$ to $p$ is reached lie on $\partial\mathrm{conv}(S)$, and one of these points is a vertex of $\mathrm{conv}(S)$.
\item If $a,b,c,d$ are the cyclically ordered vertices of a convex quadrilateral, then $\|a-c\|+\|b-d\|\geq \mathrm{max}\{\|a-b\|+\|c-d\|,\ \|a-d\|+\|b-c\|\}.$
\item Let $a,b,c$ be the vertices of a triangle. If $x\in \overline{ab}$ then $\|c-x\|\leq\max\{\|c-a\|,\|c-b\|\}$.
\end{enumerate}
\end{lemma}
\begin{proof}
(1)\; Let $q\in S$ be a point such that $\|p-q\|=\max\{\|p-x\|:x\in S\}$. Therefore, $S\subseteq B(p,\|p-q\|)$. Let $y\in \mathbb{M}^2$ be such that $q-p \perp y$. The line $q+\lambda y$ ($\lambda\in \mathbb{R}$) supports $B(p,\|p-q\|)$ at $q$, $\mathrm{conv}(S)$ and $B(p,\|p-q\|)$ are completely contained in one of the half-planes determined by the line, and any point $x$ contained in the other half-plane verifies $\|p-x\|\geq \|p-q\|$. This implies that $q$ belongs to $\partial\mathrm{conv}(S)$. If $q$ is not a vertex of $\mathrm{conv}(S)$ the segment of $\partial\mathrm{conv}(S)$ containing $q$ belongs to the line $q+\lambda y$, and the distance from $p$ to the vertices of such a segment is not smaller that $\|p-q\|$.

\noindent (2)-(3)\; See Proposition 7 and Lemma 5 in \cite{MSW}.
\end{proof}


The following statement is proved in \cite{MPSY} (Lemma 3).
\begin{lemma}\label{subset}
	Let $G=\{V,E,w:E\to \mathbb{R}^+\}$ be a weighted graph, where $V$ is the set of vertices and $E$ is the set of edges, and let $V'\subseteq V$ be a subset of $V$. Among all edges of $G$ having one endpoint in $V'$ and the other endpoint in $V\setminus V'$, let $e$ be an edge of maximum weight. Then $e$ occurs in some maximum spanning tree of $G$.
\end{lemma}

\section{The geometric structure of maximum spanning trees}\label{structure}

Since the distances among any pair of points of $S$ are distinct, every vertex of $\mathtt{FNG}(S)$ has out-degree one (every vertex is the starting point of a unique directed edge) and every (directed) path is a sequence of edges of strictly increasing length until a cycle of two edges is reached. If $x,x'$ are the vertices of such a cycle, $\overline{xx'}$ is a \textit{spine}, and this spine determines a component $\textbf{C}_{xx'}$ of $\mathtt{FNG}(S)$. A point only can be a vertex of one spine. We say that $p\in S$ belongs to the \textit{cluster} $C_x$ (respectively, cluster $C_{x'}$), if there exists a directed path from $p$ to $x$ (respectively, from $p$ to $x'$) with a even number of edges. Every point of $S$ belongs to a unique cluster, and every component of $\mathtt{FNG}(S)$ consists of exactly two clusters (each cluster is related to one vertex of the spine). Note that $x\in C_{x}$ and $x'\in C_{x'}$.

\begin{lemma}\label{onboundary} For a set $S=\{p_1,p_2,\dots,p_n\}$, the following holds:
\begin{enumerate}
    \item Only the leaves of $\mathtt{FNG}(S)$ can be internal to $\mathrm{conv}(S)$.
	\item Both endpoints of a spine of $\mathtt{FNG}(S)$ lies on $\partial\mathrm{conv}(S)$ and each pair of spines of $\mathtt{FNG}(S)$ properly intersects.
	\item The endpoints of the spines of $\mathtt{FNG}(S)$ can be labeled in cyclic order around $\partial\mathrm{conv}(S)$ as $p_1,p_2,...,p_k, q_1,q_2,...,q_k,$ where the spines are $\overline{p_iq_i}$, for $i=1,2,...,k.$
	\item All the points of a cluster that lie on $\partial\mathrm{conv}(S)$ are contiguous (i.e., if $a,\bar{a}$ belong to  a cluster and $b,\bar{b}$ belong to a different cluster, then the clockwise order around $\partial\mathrm{conv}(S)$ cannot be $ab\bar{a}\bar{b}$).
\end{enumerate}
\end{lemma}
\begin{proof}
Lemma \ref{furconvex}(1) proves (1). Every spine joints two furthest neighbours, that by Lemma \ref{furconvex}(1) must be on the $\partial\mathrm{conv}(S)$. Therefore (1) and Lemma \ref{furconvex}(2) implies (2). (3) is a consequence of (2).

Let us assume that (4) fails. Then there exist two components $\textbf{C}_{xx'}$ and $\textbf{C}_{yy'}$, and four points $a,b,\bar{a},\bar{b}\in S$ clockwise situated on $\partial\mathrm{conv}(S)$, where $a,\bar{a}\in C_x$   
and $b,\bar{b}\in C_y$.


Let us consider the sequence $(a_0:=a,a_1,...)$ starting from $a$ and where $a_{n+1}$ is the furthest neighbour of $a_n$. Such a sequence reaches $x$ in an even number of steps, and after that, $x$ and $x'$ alternate. Similarly we denote $(\bar{a}_0:=\bar{a},\bar{a}_1,...)$, $(b_0:=b,b_1,...)$, and $(\bar{b}_0:=\bar{b},\bar{b}_1,...)$ the sequences of furthest points starting from $\bar{a},b,$ and $\bar{b}$, respectively, which reach $x$, $y$, and $y$, respectively, in an even number of steps each.

Since $\overline{a_0a_1}$, $\overline{b_0b_1}$, $\overline{\bar{a}_0\bar{a}_1}$ and $\overline{\bar{b}_0\bar{b}_1}$ are  (directed) edges of $\mathtt{FNG}(S)$, $a_0a_1b_1b_0$ (in this cyclic order) is not a convex quadrilateral (Lemma \ref{furconvex}(2); see Lemma 4 in \cite{MPSY} for details), and neither $a_0a_1\bar{a}_1\bar{a}_0$, $b_0b_1\bar{a}_1\bar{a}_0$, $\bar{a}_0\bar{a}_1\bar{b}_1\bar{b}_0$, nor $b_0b_1\bar{b}_1\bar{b}_0$. Therefore the clockwise cyclic order $a_0b_0\bar{a}_0$ implies the clockwise cyclic order $a_1b_1\bar{a}_1$ on $\partial\mathrm{conv}(S)$, and the clockwise cyclic order $a_0b_0\bar{a}_0\bar{b}_0$ implies the clockwise cyclic order $a_1b_1\bar{a}_1\bar{b}_1$ on $\partial\mathrm{conv}(S)$. Repeating the argument, finally we conclude the clockwise cyclic order $a_nb_n\bar{a}_n\bar{b}_n$  on $\partial\mathrm{conv}(S)$ for any $n$.
Since $a,\bar{a}\in C_x$, for some step $n$ the sequences $a_n$ and $\bar{a}_n$ coincide in $x$, and  $b_n$ and $\bar{b}_n$ coincide in $y$. Nevertheless these coincidences and the clockwise cyclic order $a_nb_n\bar{a}_n\bar{b}_n$ are not compatible, and we get a contradiction.

\end{proof}

\begin{lemma}\label{bisector}
Let $\mathbb{M}^2$ be a normed plane and let $p,q\in\mathbb{M}^2$ be two different points. There exist a simple curve $Bi(p, q)\subseteq \mathrm{bis}(p, q)$ through the midpoint of $p$ and $q$, which is symmetric with respect to this midpoint, homeomorphic to $\mathbb{R}$, and such that for every $r\in Bi(p, q)$
the curve $Bi(p,q)$ belongs to the double cone defined by the lines $L(r,p)$ and $L(r,q)$ that contains $\frac{1}{2}(p+q)$. This curve separates the plane into two connected parts, $P(p,q)$ and $P(q,p)$, such that $\|p-y\|\leq\|q-y\|$ for $y\in P(p,q)$, and $\|q-y\|\leq\|p-y\|$ for $y\in P(q,p)$. Moreover, for any $r\in\mathrm{bis}(p, q)$ the curve $Bi(p, q)$ can be constructed to pass through $r$. If $\mathbb{M}^2$ is strictly convex, every point of $\mathrm{bis}(p, q)$ different from $r$ is an interior point of the cone.
\end{lemma}
%

\begin{proof}
The  statements are proved in 
\cite{MSW} for any strictly convex normed plane and in \cite{JS} for any normed plane, except what refers to the last sentence that we prove as follows.
	We shall use the fact that for any $x,y\in\mathbb{M}^2$ the function $\mu\in\mathbb{R}\mapsto \|x+\mu y\|$ is convex, and if $\mathbb{M}^2$ is strictly convex this function cannot be constant in any interval. Assume that there exists $t\in\mathrm{bis}(p, q) $ that is not an interior point of the cone, i.e., there exists $\lambda\in\mathbb{R}$, $0\neq\lambda\neq 1$, such that
$t=\lambda p+(1-\lambda)r$ and $\|t-p\|=\|t-q\|$. The function $\mu\in\mathbb{R}\mapsto f(\mu)=\|(1-\lambda)(r-p)+\mu(p-q)\|$ satisfies $f(0)=|1-\lambda|\|r-p\|$, $f(1)=\|t-q\|=\|t-p\|=|1-\lambda|\|r-p\|$ and $f(1-\lambda)=|1-\lambda|\|r-q\|$. Then, $f(0)=f(1)=f(1-\lambda)$, and therefore $f(\mu)$ is constant in the interval that contains the point $0$, $1$, and $1-\lambda$, which is absurd. 

\end{proof}

\begin{lemma}\label{key}Let $\textbf{C}_{xx'}$, $\textbf{C}_{yy'}$ and $\textbf{C}_{zz'}$ be three components of $\mathtt{FNG}(S)$, such that $C_x,C_y,C_z,C_{x'},C_{y'}, C_{z'}$ is the clockwise order. Then the maximum distance from a point of $\textbf{C}_{xx'}$ to $\textbf{C}_{yy'} \cup \textbf{C}_{zz'}$ is realized between a point of $\textbf{C}_{xx'}$ and a point of a non contiguous cluster (either between $C_x$ and $C_z\cup C_{y'}$ or between $C_{x'}$ and $C_y\cup C_{z'}$).
\end{lemma}
\begin{proof}
The structure of our proof is identical to that in \cite{MPSY} for the Euclidean case. For the clarity of the exposition, we present all the steps but emphasizing the arguments that are different in $\mathbb{E}^2$. Let us suppose that the lemma fails. Without loss of generality, let us assume that there exist $a'\in C_{x'}$ and $c\in C_z$ such that $\|a'-c\|$ is the maximum distance between points of $\textbf{C}_{xx'}$ and $\textbf{C}_{yy'} \cup \textbf{C}_{zz'}$. Let $\overline{bb'}$ the spine of $\textbf{C}_{yy'}$, where $b\in C_{y}$ and $b'\in C_{y'}$. Let $c'\in C_{z'}$ be the furthest neighbour of $c$ in $S$ and let $a\in C_x$ be the furthest neighbour of $a'$. We hold the following:
\begin{enumerate}
\item $\|a'-a\|>\|a'-c\|>\|a'-b\|.$
\item $\|a'-c\|>\|a-c\|, \|a'-c\|>\|b'-a\|.$
\item $\|b'-b\|>\|b'-a\|,\|b'-b\|>\|b'-c\|.$

\item \textit{The points $a,b,b',c'$ lie on $\partial\mathrm{conv}(S)$ in this cyclic order}. This follows by (1) in Lemma \ref{onboundary} and the cyclic order of  $C_x,C_y,C_z,C_{x'},C_{y'}, C_{z'}$.

\item \textit{$c$ lies on $\partial\mathrm{conv}({b,c,a',b',c'})$}. This follows by observing that $c$ is further from $a'$ than any other point in $\{b,c,a',b',c'\}$ and Lemma \ref{furconvex}(1).

\item \textit{$b,c,b',c'$ are the cyclically ordered vertices of a convex quadrilateral}. This follows by Lemma \ref{furconvex} because $b,b',c'$ is the clockwise order on $\mathrm{conv}(S)$ and $c$ is further from $c$ than any other point in $\{b,b',c'\}$.

\item \textit{$a',a,b,c$ are cyclically ordered vertices of a convex quadrilateral}. This follows because (4) and (6) guarantee that $a,b,c,b',c'$ are the cyclically ordered vertices of a convex pentagon, and $a,c,a',c$ must be the cyclically ordered vertices of a convex quadrilateral by Lemma \ref{furconvex}.

\end{enumerate}
Let $\mathrm{bis}(a, b)$ and $\mathrm{bis}(b, c)$ be the bisectors defined by the points $a,b$ and $b,c$, respectively. Let us choose a point $O\in \mathrm{bis}(a,b)\cap \mathrm{bis}(b, c)$. By Lemma \ref{bisector}, we can choose three simple curves $Bi(a,c)$, $Bi(b,c)$, and $Bi(a,b)$, each of them homeomorphic to $\mathbb{R}$, that satisfy the following conditions:
\begin{enumerate}
	\item[I] $Bi(a, c)\subseteq \mathrm{bis}(a, c)$ goes through the midpoint of $a$ and $c$, is symmetric with respect to this midpoint, and belongs to the double cone with apex $O$ and through $a$ and $c$. This curve separates the plane into two connected parts $P(c, a)$ and $P(a,c)$ such that whenever $y\in P(c, a)$, then $\|c-y\|\leq \|a-y\|$, and whenever $y\in P(a, c)$, then $\|a-y\|\leq \|c-y\|.$
	\item[II] $Bi(b, c)\subseteq \mathrm{bis}(b, c)$ goes through the midpoint of $b$ and $c$, is symmetric with respect to this midpoint, and belongs to the double cone with apex $O$ and through $b$ and $c$. This curve separates the plane into two connected parts $P(c, b)$ and $P(b,c)$ such that whenever $y\in P(c, b)$, then $\|c-y\|\leq \|b-y\|$, and whenever $y\in P(b, c)$, then $\|b-y\|\leq \|c-y\|.$
	\item[III] $Bi(a, b)\subseteq \mathrm{bis}(a, b)$ goes through the midpoint of $a$ and $b$, is symmetric with respect to this midpoint, and belongs to the double cone with apex $O$ and through $a$ and $b$. This curve separates the plane into two connected parts $P(b, a)$ and $P(a,b)$ such that whenever $y\in P(b,a)$, then $\|b-y\|\leq \|a-y\|$, and whenever $y\in P(a, b)$, then $\|a-y\|\leq \|b-y\|.$
\end{enumerate}
The following holds:
\begin{enumerate}
	\item[(A)] \textit{$a'$ lies in the triangle with vertices $\{a,O,c\}.$} This follows from the fact $a'\in P(c,a)\cap P(b,c)\cap P(b,a)$ and the inclusion of  $Bi(a,c)$ inside the double cone with apex $O$ and through $a$ and $c$, the inclusion of $Bi(b,c)$ inside the double cone with apex $O$ and through $b$ and $c$, and the inclusion of $Bi(a,b)$ inside the double cone with apex $O$ and through $a$ and $b$.
	
	\item[(B)] \textit{$b'$ lies in the "wedge" that is the intersection of right half-plane of the directed half-line from $a$ to $O$ and the left-plane of the directed half-line from $c$ to $O$}. This follows from the fact $b'\in P(a,b)\cap P(c,b)$ and the inclusion of  $Bi(a,b)$  inside the double cone with apex $O$ and through $a$ and $b$, and the inclusion of $Bi(b,c)$ inside the double cone with apex $O$ and through $b$ and $c$.
\end{enumerate}
(A) and (B) imply that $a'\in \mathrm{conv}(\{a,c,b'\})$, therefore $\|a-a'\|\leq \mathrm{max}\{\|a-c\|,\|a-b'\|\}$ by Lemma \ref{furconvex}(3), and we get a contradiction with (1) and (2).
\end{proof}

%
%
%
%

We prove now that all the edges of $\mathtt{MXST}(S)$ that are not in $\mathtt{FNG}(S)$ join adjacent components of $\mathtt{FNG}(S)$, in the cyclic ordering of components.
\begin{lemma}\label{ordercomponents}
	Let $\textbf{C}_{x_{i-1}x_{i-1}'}$, $\textbf{C}_{x_{i}x'_{i}}$ and $\textbf{C}_{x_{i+1}x'_{i+1}}$ be three components such that the clusters $C_{x_{i-1}}$ and $C_{x_{i+1}}$ are adjacent to $C_{x_i}$ and the cyclic ordering is $C_{x_{i-1}},C_{x_i},C_{x_{i+1}}$. If there is an edge $e\in \mathtt{MXST}(S)\setminus \mathtt{FNG}(S)$ with one endpoint in $C_{x_i}$, then the other endpoint of $e$ is either at $C_{x'_{i-1}}$ or at $C_{x'_{i+1}}$.
\end{lemma}
\begin{proof}
Note that if $e\in \mathtt{MXST}(S)\setminus \mathtt{FNG}(S)$ then $e$ joins points from two different components. Suppose that $e$ does not satisfy the condition of the lemma. Without of loss of generality we assume that $e$ has an endpoint in $C_{y}$ such that the cyclic ordering is $C_{x_{i-1}},C_{x_i},C_{x_{i+1}}, C_y$.  By Lemma \ref{key}, the length of $e$ is not the maximum distance between points of $\textbf{C}_{x_{i}x'_{i}}$ and $\textbf{C}_{x_{i-1}x'_{i-1}}\cup \textbf{C}_{yy'}$. We obtain a spanning tree larger than $\mathtt{MXST}(S)$ replacing $e$ by the edge joining the pair of points where this maximum distance is reached, which is a contradiction.
\end{proof}

\begin{lemma}\label{main}
	If $e$ is an edge of $\mathtt{MXST}(S)$, then at least one endpoint of $e$ lies on $\partial\mathrm{conv}(S)$.
\end{lemma}

\begin{proof}

Let us assume that $e=\overline{xy}$ is an edge of $\mathtt{MXST}(S)$ such that $x,y$ are not on the boundary of $\mathrm{conv}(S)$. Let $A$ and $B$ be the two components of $\mathtt{MXST}(S)$ obtained when $e$ is deleted, where $x\in A$ and $y\in B$. We prove in the following that there exist $u\in A$ and $v\in B$ such that $\|u-v\|>\|x-y\|$ and at least one of these points is on $\partial\mathrm{conv}(S)$. As a consequence, a larger spanning tree could be obtained replacing $e$ by $\overline{uv}$, which is a contradiction.

\bigskip

\begin{figure}[h]
	\begin{center}
		\scalebox{0.8}{\includegraphics{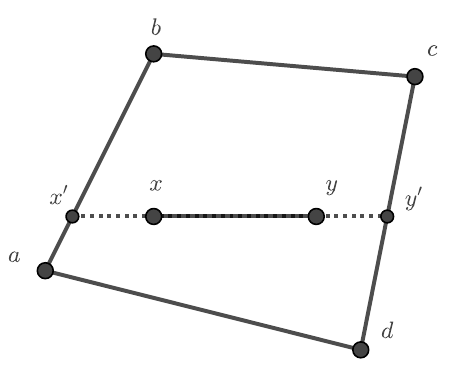}}\\
		\caption{\label{eboundary}$\overline{xy}\notin \mathtt{MXST}(S)$ if $x,y$ are not on $\partial\mathrm{conv}(S)$.}
	\end{center}
\end{figure}

\bigskip

Let $x'$, $y'$ be the points on the boundary of $\mathrm{conv}(S)$ such that $x'-y=\lambda (x-y)$ with $\lambda >1$ and $y'-x=\beta (y-x)$ with $\beta >1$.
Let $\{a,b,c,d\}$ be the cyclically ordered  vertices on the boundary of $\mathrm{conv}(S)$ such that $x'\in \overline{ab}$, and $y'\in \overline{cd}$ (see Figure \ref{eboundary}).

Next we shall consider all the possible situations depending on which component either $A$ or $B$ the points $\{a,b,c,d\}$ are located. This will give rise to sixteen cases, that we identify with the following notation: e.g., if $a\in A$, $b\in B$, $c\in B$ and $d\in B$, this is the case $ABBB$; the cases in which we only assume that $a,b\in A$ will be denoted by $AAXX$; and so on.
Note that by interchanging the points $\{a,b,c,d\}$ and the sets $A, B$ some cases are equivalent (see Table \ref{tab:table}), which reduces the study to the cases $AAXX$, $ABAA$, $ABAB$, $ABBA$ and $BBAA$.
{\footnotesize
\begin{table}
\renewcommand{\arraystretch}{1.1}
\begin{tabular}{c|c|c|c|c|c|}
\cline{2-6}
& \multicolumn{5}{|c|}{Equivalent to}\\\hline
\multicolumn{1}{|c|}{Case} & AAXX      & ABAA      & ABAB      & ABBA      & BBAA      \\\hline
\multicolumn{1}{|l|}{AAAA} & $\bullet$ &           &           &           &           \\\hline
\multicolumn{1}{|l|}{AAAB} & $\bullet$ &           &           &           &           \\\hline
\multicolumn{1}{|l|}{AABA} & $\bullet$ &           &           &           &           \\\hline
\multicolumn{1}{|l|}{AABB} & $\bullet$ &           &           &           &           \\\hline
\multicolumn{1}{|l|}{ABAA} &           & $\circ$   &           &           &           \\\hline
\multicolumn{1}{|l|}{ABAB} &           &           & $\circ$   &           &           \\\hline
\multicolumn{1}{|l|}{ABBA} &           &           &           & $\circ$   &           \\\hline
\multicolumn{1}{|l|}{ABBB} & $\bullet$ &           &           &           &           \\\hline
\multicolumn{1}{|l|}{BAAA} &           & $\bullet$ &           &           &           \\\hline
\multicolumn{1}{|l|}{BAAB} &           &           &           & $\bullet$ &           \\\hline
\multicolumn{1}{|l|}{BABA} &           &           & $\bullet$ &           &           \\\hline
\multicolumn{1}{|l|}{BABB} & $\bullet$ &           &           &           &           \\\hline
\multicolumn{1}{|l|}{BBAA} &           &           &           &           & $\circ$   \\\hline
\multicolumn{1}{|l|}{BBAB} &           & $\bullet$ &           &           &           \\\hline
\multicolumn{1}{|l|}{BBBA} &           & $\bullet$ &           &           &           \\\hline
\multicolumn{1}{|l|}{BBBB} & $\bullet$ &           &           &           &           \\\hline
\end{tabular}
\bigskip
\caption{Cases in Lemma \ref{main}}
\label{tab:table}
\end{table}
}

To simplify the counts, without loss of generality, we can assume that $x=(0,0)$, $y=(1,0)$, and $\|y\|=1$; as well that $a=(a_1,a_2), b=(0,b_2), c=(c_1,c_2), d=(d_1,d_2),$ where
$$
a_1<0,\quad a_2\leq0,\quad b_2>0,\quad c_1\geq0,\quad c_2\geq0,\quad d_1\in\mathbb{R},\quad d_2\leq0.
$$
Also note that if $d_1<1$ then $c_1>1$. Moreover, since $x$ and $y$ are not on $\partial\mathrm{conv}(S)$, $a_2$ and $d_2$ cannot be simultaneously zero.



Since $y'=\big(\frac{c_2d_1-c_1d_2}{c_2-d_2},0\big)=\beta y$, with $\beta>1$, we have that
\begin{equation}\label{sign_c_d}
h:=c_2d_1-c_1d_2+d_2-c_2> 0.\tag{$*$}
\end{equation}


\noindent Cases $AAXX$: By Lemma 2.1(3) we have $\|x-y\|\leq\max\{\|a-y\|,\|b-y\|\}$. Thus we take $v=y$ and either $u=a$ or $u=b$.

%

\medskip

In the remaining cases we will suppose that $\|u-v\|<\|x-y\|$ for every $u\in A$ and $v\in B$, $\{u,v\}\neq\{x,y\}$, and we will get a contradiction. For that, in any of the different cases we will show that
$$
y=\alpha_1 x+\alpha_2(u-v)+\alpha_3(\bar u-\bar v),
$$
with $\alpha_1,\alpha_2,\alpha_3\in[0,1]$, $\alpha_1+\alpha_2+\alpha_3=1$, being $\{u,v\}$ and $\{\bar u,\bar v\}$ either in $A\times B$ or in $B\times A$. This means that $y$ is inside the triangle of vertices $x$, $u-v$, $\bar u-\bar v$, and then by Lemma 2.1(3), $\|x-y\|\leq\max\{\|u-v\|,\|\bar u-\bar v\|\}$, contradicting the initial supposition.

In all the cases we will write $\alpha_i=\frac{s_i}{s_i+t_i}$, $i=1,2,3$. Because of the sign of the parameters involved, we have that $s_i,t_i\geq 0$, which implies $\alpha_i\in[0,1]$. Moreover, $s_1+t_1=s_2+t_2=s_3+t_3>0$, and $\alpha_1+\alpha_2+\alpha_3=1$. In some cases we will use $h$ as it is defined in (\ref{sign_c_d}).

\medskip

\noindent Case $ABAA$:

\noindent (1) Assume that $d_1\geq1$. We have $y=\alpha_1 x+\alpha_2 (y-a)+\alpha_3 (d-b)$ with
\begin{align*}
\alpha_1 & =\frac{[-a_2(d_1-1)-a_1(b_2-d_2)]}{[-a_2(d_1-1)-a_1(b_2-d_2)]+[b_2-a_2-d_2]},\\[.5ex]
\alpha_2 & =\frac{[b_2-d_2]}{[b_2-d_2]+[a_1d_2-a_2d_1-a_1b_2]},\\[.5ex]
\alpha_3 & =\frac{[-a_2]}{[-a_2]+[-a_2(d_1-1)-(a_1-1)(b_2-d_2)]}.
\end{align*}

\medskip

\noindent (2) Assume that $d_1<1$. Then $c_1>1$.

\smallskip

(2.1)\; Assume that $c_2<b_2$. We have $y=\alpha_1x+\alpha_2(y-a)+\alpha_3(c-b)$ with
\begin{align*}
\alpha_1 & =\frac{[-a_2(c_1-1)-a_1(b_2-c_2)]}{[-a_2(c_1-1)-a_1(b_2-c_2)]+[(b_2-c_2)-a_2]},\\[.5ex]
\alpha_2 & =\frac{[b_2-c_2]}{[b_2-c_2]+[-a_1(b_2-c_2)-a_2c_1]},\\[.5ex]
\alpha_3 & =\frac{[-a_2]}{[-a_2]+[-a_2(c_1-1)-(a_1-1)(b_2-c_2)]}.
\end{align*}

\smallskip

(2.2)\; Assume that $c_2\geq b_2$. We have $y=\alpha_1x+\alpha_2(d-b)+\alpha_3(c-b)$ with
\begin{align*}
\alpha_1 & =\frac{[h+b_2(c_1-d_1)]}{[h+b_2(c_1-d_1)]+[c_2-d_2]},\\[.5ex]
\alpha_2 & =\frac{[c_2-b_2]}{[c_2-b_2]+[h+b_2(c_1+1-d_1)-d_2]},\\[.5ex]
\alpha_3 & =\frac{[b_2-d_2]}{[b_2-d_2]+[h+(c_2-b_2)+b_2(c_1-d_1)]}.
\end{align*}

\medskip

\noindent Case $ABAB$:

\noindent (1) Assume that $d_1\geq1$. We have $y=\alpha_1 x+\alpha_2 (y-a)+\alpha_3 (d-x)$ with
\begin{align*}
\alpha_1 & =\frac{[a_1d_2-a_2(d_1-1)]}{[a_1d_2-a_2(d_1-1)]+[-a_2-d_2]},\\[.5ex]
\alpha_2 & =\frac{[-d_2]}{[-d_2]+[d_2a_1-d_1a_2]},\\[.5ex]
\alpha_3 & =\frac{[-a_2]}{[-a_2]+[-a_2(d_1-1)-d_2(1-a_1)]}.
\end{align*}

\medskip

\noindent (2) Assume that $d_1<1$. Then $c_1>1$.

\smallskip

(2.1)\; Assume that $c_2< b_2$. We have $y=\alpha_1x+\alpha_2(y-a)+\alpha_3(c-b)$ with
\begin{align*}
\alpha_1 & =\frac{[-a_1(b_2-c_2)-a_2(c_1-1)]}{[-a_1(b_2-c_2)-a_2(c_1-1)]+[(b_2-c_2)-a_2]},\\[.5ex]
\alpha_2 & =\frac{[b_2-c_2]}{[b_2-c_2]+[-a_1(b_2-c_2)-a_2c_1]},\\[.5ex]
\alpha_3 & =\frac{[-a_2]}{[-a_2]+[(b_2-c_2)(1-a_1)-a_2(c_1-1)]}.
\end{align*}

\smallskip

(2.2)\; Assume that $c_2\geq b_2$. Then $y=\alpha_1 x+\alpha_2 d+\alpha_3(c-b)$ with
\begin{align*}
\alpha_1 & =\frac{[h+b_2(1-d_1)]}{[h+b_2(1-d_1)]+[(c_2-b_2)-d_2]},\\[.5ex]
\alpha_2 & =\frac{[c_2-b_2]}{[c_2-b_2]+[h-d_2+b_2(1-d_1)]},\\[.5ex]
\alpha_3 & =\frac{[-d_2]}{[-d_2]+[h+(c_2-b_2)+b_2(1-d_1)]}.\\
\end{align*}

\smallskip

\noindent Case $ABBA$:

\noindent (1) Assume that $d_1\geq1$. We have $y=\alpha_1 x+\alpha_2 (y-a)+\alpha_3 (d-b)$ with
\begin{align*}
\alpha_1 & =\frac{[-a_1(b_2-d_2)-a_2(d_1-1)]}{[-a_1(b_2-d_2)-a_2(d_1-1)]+[b_2-a_2-d_2]},\\[.5ex]
\alpha_2 & =\frac{[b_2-d_2]}{[b_2-d_2]+[-a_1(b_2-d_2)-a_2d_1]},\\[.5ex]
\alpha_3 & =\frac{[-a_2]}{[-a_2]+[(b_2-d_2)(1-a_1)-a_2(d_1-1)]}.
\end{align*}

\smallskip

\noindent (2) Assume that $d_1<1$. Then $c_1>1$. We have $y=\alpha_1 x+\alpha_2 (c-x)+\alpha_3 (d-b)$ with
\begin{align*}
\alpha_1 & =\frac{[h+b_2(c_1-1)]}{[h+b_2(c_1-1)]+[b_2+c_2-d_2]},\\[.5ex]
\alpha_2 & =\frac{[b_2-d_2]}{[b_2-d_2]+[h+b_2(c_1-1)+c_2]},\\[.5ex]
\alpha_3 & =\frac{[c_2]}{[c_2]+[h+b_2c_1-d_2]}.
\end{align*}

\smallskip

\noindent Case $BBAA$:

\noindent (1) Assume that $d_1\geq1$.

\smallskip

(1.1)\; Assume that $c_1\geq1$. We have $y=\alpha_1x+\alpha_2(c-a)+\alpha_3(d-b)$ with
\begin{align*}
\alpha_1 & =\frac{[h+a_1d_2-a_2(d_1-1)+b_2(c_1-a_1-1)]}{[h+a_1d_2-a_2(d_1-1)+b_2(c_1-a_1-1)]+[b_2+c_2-a_2-d_2]},\\[.5ex]
\alpha_2 & =\frac{[b_2-d_2]}{[b_2-d_2]+[h+(a_1d_2-a_2d_1)+b_2(c_1-1-a_1)+c_2]},\\[.5ex]
\alpha_3 & =\frac{[c_2-a_2]}{[c_2-a_2]+[h-a_2(d_1-1)-d_2(1-a_1)+b_2(c_1-a_1)]}.
\end{align*}

\smallskip

(1.2)\; Assume that $c_1<1$ and $a_2\geq d_2$. We have $y=\alpha_1x+\alpha_2(c-a)+\alpha_3(d-a)$ with
\begin{align*}
\alpha_1 & =\frac{[h-a_2(d_1-c_1)+a_1d_2-a_1c_2]}{[h-a_2(d_1-c_1)+a_1d_2-a_1c_2]+[c_2-d_2]},\\[.5ex]
\alpha_2 & =\frac{[a_2-d_2]}{[a_2-d_2]+[h+(a_1d_2-a_2d_1)+c_2(1-a_1)-a_2(1-c_1)]},\\[.5ex]
\alpha_3 & =\frac{[c_2-a_2]}{[c_2-a_2]+[h+(a_2-d_2)-a_2(d_1-c_1)-a_1(c_2-d_2)]}.
\end{align*}

\smallskip

(1.3)\; Assume that $c_1<1$ and $a_2<d_2$. We have $y=\alpha_1x+\alpha_2(d-b)+\alpha_3(d-a)$ with
\begin{align*}
\alpha_1 & =\frac{[b_2(d_1-1-a_1)-a_2(d_1-1)+a_1d_2]}{[b_2(d_1-1-a_1)-a_2(d_1-1)+a_1d_2]+[b_2-a_2]},\\[.5ex]
\alpha_2 & =\frac{[d_2-a_2]}{[d_2-a_2]+[b_2(d_1-a_1)-a_2(d_1-1)-d_2(1-a_1)]},\\[.5ex]
\alpha_3 & =\frac{[b_2-d_2]}{[b_2-d_2]+[(d_2-a_2)+(d_1-1)(b_2-a_2)-a_1(b_2-d_2)]}.
\end{align*}

\smallskip

\noindent (2) Assume that $d_1<1$. Then $c_1>1$.

(2.1)\; Assume that $c_2\geq b_2$. We have $y=\alpha_1x+\alpha_2(c-b)+\alpha_3(d-b)$ with
\begin{align*}
\alpha_1 & =\frac{[h+b_2(c_1-d_1)]}{[h+b_2(c_1-d_1)]+[c_2-d_2]},\\[.5ex]
\alpha_2 & =\frac{[b_2-d_2]}{[b_2-d_2]+[h+(c_2-b_2)+b_2(c_1-d_1)]},\\[.5ex]
\alpha_3 & =\frac{[c_2-b_2]}{[c_2-b_2]+[h+b_2(1+c_1-d_1)-d_2]}.
\end{align*}

\smallskip

(2.2)\; Assume that $c_2<b_2$. We have $y=\alpha_1x+\alpha_2(c-b)+\alpha_3(c-a)$ with
\begin{align*}
\alpha_1 & =\frac{[(b_2-a_2)(c_1-1)-a_1(b_2-c_2)]}{[(b_2-a_2)(c_1-1)-a_1(b_2-c_2)]+[b_2-a_2]},\\[.5ex]
\alpha_2 & =\frac{[c_2-a_2]}{[c_2-a_2]+[(b_2-c_2)(1-a_1)+(b_2-a_2)(c_1-1)]},\\[.5ex]
\alpha_3 & =\frac{[b_2-c_2]}{[b_2-c_2]+[b_2(c_1-1)-a_1(b_2-c_2)-a_2c_1+c_2]}.
\end{align*}

\end{proof}

\section{The Algorithm}\label{algorihtm}

\begin{figure}[h]
	\begin{center}
        \scalebox{1}{\includegraphics{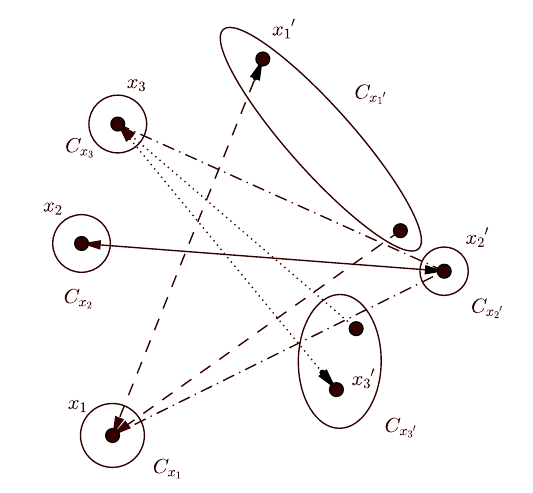}}\hspace{0.4cm}
\\\caption{\label{constructionmxstree}The construction of $\mathtt{MXST}(S)$.}
	\end{center}
\end{figure}

Given a set of points $S$ in $\mathbb{E}^2$ with distinct distances among any pair of points, the algorithm by Monma et al. for computing the $\mathtt{MXT}(S)$ runs as follows (see Figure \ref{constructionmxstree}).

\bigskip


\begin{algorithm}

 \label{maximum}
\begin{enumerate}
\item For each point $p\in S$ compute its furthest neighbour in $S$.
\item Construct the directed graph $\mathtt{FNG}(S)$, whose nodes are the points of $S$ and
whose (directed) edges connect points of $S$ to their respective furthest
neighbours.
\item Let $\textbf{C}_{x_1x_1'}, \textbf{C}_{x_2x_2'}, \dots, \textbf{C}_{x_kx_k'}$, be the components of $\mathtt{FNG}(S)$, where the clusters $C_{x_1}, C_{x_2},\dots, C_{x_k}, C_{x_1'},C_{x_2'},\dots, C_{x_k'}$ are cyclically ordered around $\mathrm{conv}(S)$.
\item For $i= 1 , 2 , . . . , k$, compute the edge (and its extreme points) with the following length (if $i$=1, then $i-1:=k$):
$$e_{i}=\max
\left\{\begin{array}{c}
\max\{\|x-y\|\;/\; x \in C_{x_i},\; y\in C_{x'_{i-1}}\cap \partial\mathrm{conv}(S)\},\\
\max\{\|x-y\|\;/\; x \in C_{x_{i-1}},\; y\in C_{x'_{i}}\cap \partial\mathrm{conv}(S)\},\\
\max\{\|x-y\|\;/\; x \in C_{x'_i},\; y\in C_{x_{i-1}}\cap \partial\mathrm{conv}(S)\},\\
\max\{\|x-y\|\;/\; x \in C_{x'_{i-1}},\; y\in C_{x_{i}}\cap \partial\mathrm{conv}(S)\},\\
\end{array}\right.
$$
Discard the edge with the smallest length. Let $E$ be the set of  remaining $k-1$ edges.
\item Output $\mathtt{MXT}(S)$ as the union of the edges of $\mathtt{FNG}(S)$ and $E$.
\end{enumerate}
\end{algorithm}

\bigskip

\begin{thm} Let $\mathbb{M}^2$ be a normed plane. For a given set $S$ of $n$ points with distinct distances among pair of points, the $\mathtt{MXST}(S)$ can be computed using Algorithm \ref{maximum}.
\end{thm}
\begin{proof}
Algorithm \ref{maximum} locates  all edges of $\mathtt{FNG}(S)$ in Step 2. These edges are in $\mathtt{MXST}(S)$ by Lemma \ref{subset}. Algorithm \ref{maximum} locates the remaining edges of $\mathtt{MXST}(S)$ in Step 4, which must connect points of two different components of $\mathtt{FNG}(S)$, based on Lemma \ref{ordercomponents} and on Lemma \ref{main} .
\end{proof}

The \textit{Voronoi diagram} and the \textit{furthest Voronoi diagram} are some of most frequently investigated structures in Computational Geometry, and the second one is useful in order to locate the furthest neighbour of a point of $S$. Klein \cite{Klein} introduces the \textit{abstract Voronoi diagram} as a unifying approach to Voronoi diagrams based on the notion of bisecting curves instead of the concept of distance. In the abstract Voronoi diagram framework the bisecting curves should be homeomorphic to a line, and each pair of them may intersect in only a finite number of points. For each pair $p$ and $q$ of sites (points, lines, etc.) the existence of a bisector $Bi(p,q)$ dividing the plane into a $p$-region $P(p,q)$ and a $q$-region $P(q,p)$ is postulated. 
Given a set $S$ of $n$ sites, the \textit{nearest abstract Voronoi} region of $p\in S$ related to $S$ is the set $\cap_{q\in S,q\neq p} P(p,q)$, namely, the locus of points closer to $p$ than to other point of $S$.
The abstract Voronoi diagram can be represented as a planar graph in a natural way (\cite{Klein}) and can be computed in expected time $O(n \log n)$ and expected space $O(n)$ (\cite{KMM}).
Mehlhorn et al. \cite{Mehlhorn} study the so-called \textit{furthest abstract Voronoi diagrams}. The furthest abstract Voronoi region of $p\in S$ related to $S$ is the locus of points further away from $p$ than from other point of $S$. They show
that furthest site abstract Voronoi diagrams are trees, have linear size, and that, given a set of $n$ sites, the  diagram can be computed by a randomized algorithm in expected time $O(n \log n)$.
Furthest abstract Voronoi diagrams encompass a large number of specific diagrams, e.g., diagrams for point sites under
any $L_p$-norm ($1 < p < \infty$).

\begin{corollary}\label{general}
	Let $\mathbb{M}^2$ be a normed plane and $S$ be a set of $n$ points with distinct distances among pair of points. The $\mathtt{MXST}(S)$ can be computed:
	\begin{enumerate}
		\item[(a)] By a randomized algorithm in expected time $O(n \log n)$ if any furthest Voronoi diagram for a set of $n$ points is a furthest abstract Voronoi diagram.
		\item[(b)] By a deterministic algorithm in expected time $O(n \log n)$ if any furthest Voronoi diagram for a set of $n$ points is a $n$-vertex planar subdivision that can be computed in $O(n\log n)$ time by a subroutine.
	\end{enumerate}
\end{corollary}

\begin{proof}
Let $S$ be a set of $n$ points in $\mathbb{M}^2$. Algorithm \ref{maximum} is suitable and the time complexity depends on the construction of the furthest neighbour of a point in a set of points. Let us consider $S_1:=S\cap\partial\mathrm{conv}(S)$. If $h$ is the number of points of $S_1$, we can compute $S_1$ in time $O(n\log h)$.

For (a), the furthest site abstract Voronoi diagram of $S_1$ can be computed by a randomized algorithm in expected time $O(h \log h)$ (\cite{Mehlhorn});
the location of any point in the diagram can be effected in  $\log h$ time after $O(h\log h)$ preprocessing time (\cite{P-S});  the construction of the $\mathtt{FNG}(S)$ therefore takes $O(n\log h)$. Each of the steps (2) and (3) requires $O(n)$ time. Step (4) takes $O(n_i\log h_i)$ time for every $i\in \{1,2,\dots,k\}$, where $n_i$ is the total number of points in $C_{x_i}\cup C_{x_{i-1}}\cup C_{x'_i}\cup C_{x_{i-1}'}$ and $h_i$ is the number of these points
lying on $\partial\mathrm{conv}(S)$. Thus Step 4 takes $O(n\log h)$ time. The total cost is $O(n\log n)$ time.
%

For (b), compute the furthest point Voronoi diagrams using the available subroutine instead of the randomized algorithm of Mehlhorn et al. for the furthest abstract Voronoi diagram.
\end{proof}

\section{Some applications: the 2-clustering problem of minimizing the maximum diameter.}\label{applications}
Given a set of points $S$ in a metric space, the \textit{$k$-clustering problem of minimizing the maximum diameter} asks how to split $S$ into $k$ sets minimizing the maximum diameter. The following algorithm computes the problem in any normed plane (\cite{Martin-Ya}) for $k=2$.

 \begin{algorithm}

 \label{Avis}Given a set $S$ of $n$ points in the plane:
	\begin{enumerate}
	\item Construct a \textit{maximum spanning tree} 
  for $S$. 
Sort the distances $d_i$ between the points of the maximum spanning tree of $S$ in increasing order.
	\item Locate the minimum $d_i$ that admits a \textit{stabbing line}\footnote{A \textit{stabbing line} for a set of segments is a line that intersects every segment of the set.}
for the maximum spanning tree by a binary search. Use the graph $(S,E_{d_i})$, where $E_{d_i}$ is the set of edges at a distance greater than $d_i$, and use the algorithm by Edelsbrunner et al. \hspace{-0.3cm} (\cite{E-M-P-R-W-W}) in order to find a stabbing line for $E_{d_i}$  as a subroutine.
	\end{enumerate}
\end{algorithm}

\begin{thm}
Let $\mathbb{M}^2$ be a normed plane and $S$ be a set of $n$ points with distinct distances among pair of points. If any furthest Voronoi diagram for a set of $n$ points is a $n$-vertex planar subdivision that can be computed in $O(n\log n)$ time by a subroutine, then the $k$-clustering problem of minimizing the maximum diameter of $S$ can be solved
 \begin{enumerate}
 \item in $O(n \log n)$ time for $k=2$,
\item  in $O(n^2 \log^2 n)$ time for $k=3$.
\end{enumerate}
\end{thm}
\begin{proof}
Algorithm \ref{Avis} works correctly for $k=2$ by Theorem 2 in \cite{Asano} and Theorem 2.9 in \cite{Martin-Ya}. The time cost of step (1) is $O(n\log n)$ by Theorem \ref{general}; finding the stabbing line takes $O(m \log m)$ times, where $m$ is the number of edges  (\cite{E-M-P-R-W-W}); the partition into two subsets in step (3) can be obtained in $O(n)$ time easily once the stabbing line is known. Therefore, we hold (1).

Algorithm 3.3 in \cite{Martin-Ya} (see also \cite{Hagauer-Rote}) computes the solution for $k=3$ and the complexity depends on the complexity of case $k=2$. The proof of (2) is exactly the proof of Theorem 3.10 in \cite{Martin-Ya} but considering that the $2$-clustering problem can be computed in $O(n\log n)$ by (1).
\end{proof}

\section{Slight movements to get different distances}\label{slightmovements}

We present in Proposition \ref{movingpoints} a descriptive step-by-step procedure for moving \textit{slightly} the points in order to break the coincidences of distances and be able to apply Algorithm \ref{maximum}, although maybe the existence of such movement could be proved by a much simpler argument. At the end of this section we justify that a similar result for non strictly convex normed planes is not true.
\smallskip

A thorough view of the proof of Lemma \ref{bisector} justifies the following.

\begin{lemma}\label{bisectorcone2}
Let $\mathbb{M}^2$ be a strictly convex normed plane. If $y\in Bi(p, q)$ for any pair of distinct points $p$ and $q$, then $\|y+\epsilon (y-p)-q\|<\|y+\epsilon (y-p)-p\|$ for any $\epsilon>0$.
\end{lemma}	


\begin{proposition}\label{movingpoints}
	Let $\mathbb{M}^2$ be a strictly convex normed plane and $S$ be a set of $n$ points. For any $\epsilon>0$ there is a point set $S'$ in a bijective correspondence with $S$ such that:
	\begin{enumerate}
		\item Each point in $S'$ is within distance $\epsilon$ of its corresponding point in $S$.
		\item The distances among points of $S'$ are all distinct.
	\end{enumerate}
\end{proposition}
\begin{proof} Let $\epsilon>0$. Let us consider the different values $\{d_1,d_2,...,d_t\}$ of the distances among the points and
 $$m=\min\{|d_r-d_s|>0: s,r\in \{1,2,...,t\}\}.$$ 
 Let $\{p_1,p_2,...,p_k, p_{k+1}, ...,p_n\}$ be the points of $S$ such that $\{p_1,p_2,...,p_k\}$ are the  vertices of  $\mathrm{conv}(S)$  ordered in counterclockwise direction. All the points in $S\cap \partial\mathrm{conv}(S)$ are considered vertices as well for our purpose. A side of $\mathrm{conv}(S)$ is a segment jointing two adjacent vertices.

Step 1. Objective: \textit{To discard that the sides of  $\mathrm{conv}(S)$ are involved in any equality of distances.}

If $p_k$ is involved in any equality of distances, take $\epsilon_{k,1}<\min\{m,\epsilon\}$
and replace $p_k$ by
$$p'_k:=p_k+\frac{\epsilon_{k,1}}{2n}\left(\frac{p_k-p_{k-1}}{\|p_k-p_{k-1}\|}\right).$$
Since $\|p_{k-1}-p'_k\|=\|p_{k-1}-p_k\|+\frac{\epsilon_{k,1}}{2n}$,  the distance $\|p_{k-1}-p'_k\|$ is not involved in any equality of distances by Lemma \ref{bisectorcone2} and the election of $\epsilon_{k,1}$. Besides, if $\|p_i-p_k\|=\|p_j-p_k\|$ for some $p_i,p_j\in S$ ($i\neq j$), the movement implies that $\|p_i-p'_k\|\neq \|p_j-p'_k\|$  by Lemma \ref{bisector} and Lemma \ref{bisectorcone2}. For simplicity, let us denote $p_k:=p'_k$. Therefore, after this movement:
\begin{itemize}
	\item The distance $\|p_{k-1}-p_k\|$ is not involved in any equality of distances.
	\item $\|p_i-p_k\|\neq\|p_j-p_k\|$ $\forall p_i,p_j\in S, i\neq j$.
\end{itemize}

Reevaluate $m$ and repeat the above argument: if $p_{k-1}$ is involved in any equality of distances, take
$\epsilon_{k-1,1}<\min\{m,\epsilon\}$
and replace $p_{k-1}$ by
$$p'_{k-1}:=p_{k-1}+\frac{\epsilon_{k-1,1}}{2n}\left(\frac{p_{k-1}-p_{k-2}}{\|p_{k-1}-p_{k-2}\|}\right).$$
Denoting $p_{k-1}:=p'_{k-1}$ and using again Lemma \ref{bisector} and Lemma \ref{bisectorcone2}, we obtain the following:
\begin{itemize}
	\item The distance $\|p_{k-2}-p_{k-1}\|$ is not involved in any equality of distances.
	\item $\|p_i-p_{k-1}\|\neq\|p_j-p_{k-1}\|$ $\forall p_i,p_j\in S, i\neq j$.
\end{itemize}
After applying similar movements to the rest of vertices of $\mathrm{conv}(S)$ as many times as necessary, we get that:
\begin{itemize}
	\item The equalities of distances do not affect to any side of $\mathrm{conv}(S)$.
	\item $\|p_{i}-p_h\|\neq\|p_{j}-p_h\|$ $\forall i,j\in \{1,2,...,n\}, i\neq j$, and  $\forall h\in\{1,...,k\}$.
\end{itemize}

\bigskip

Note that some distances like $\|p_{k-2}-p_k\|$ may not have changed after this Step 1 (even when $p_k$ or $p_{k-2}$ are moved) and $\|p_{k-2}-p_k\|$ could be involved in any equality of distances.

\bigskip

Step 2. Objective: \textit{to discard that $\|p_{h-2}-p_h\|$ is involved in any equality of distances for all $h\in\{1,...,k\}$ (denote $p_{0}:=p_k$ and $p_{-1}:=p_{k-1}$).}

Reevaluate $m$. If $p_{k-2}-p_k$ is involved in any equality of distances, take
$\epsilon_{k,2}<\min\{m,\epsilon\}$ 
and replace $p_{k}$ by
$$p'_{k}:=p_{k}+\frac{\epsilon_{k,2}}{2n}\left(\frac{p_{k}-p_{k-2}}{\|p_{k}-p_{k-2}\|}\right).$$
This movement implies that $\|p'_{k}-p_{k-2}\|$ is not involved in any equality of distance  by Lemma \ref{bisector} and Lemma \ref{bisectorcone2}. Let us denote $p_{k}:=p'_{k}$.

Applying similar movements to each vertex of  $\mathrm{conv}(S)$ as many times as necessary we achieve the objective of Step 2.

\bigskip

Step 3. Objective: \textit{to discard that $\|p_{h-3}-p_h\|$ is involved in any equality of distances for all $h\in\{1,...,k\}$ (denote $p_{0}:=p_k$, $p_{-1}:=p_{k-1}$, and $p_{-2}:=p_{k-2}$).}

Reevaluate $m$. If $p_k-p_{k-3}$ is involved in any equality of distances, take
$\epsilon_{k,3}<\min\{m,\epsilon\}$ 
and replace $p_{k}$ by $p'_{k}:=p_{k}+\frac{\epsilon_{k,3}}{2n}\left(\frac{p_{k}-p_{k-3}}{\|p_{k}-p_{k-3}\|}\right).$ This movement implies that $\|p'_{k}-p_{k-3}\|$ is not involved in any equality of distance  by Lemma \ref{bisector} and Lemma \ref{bisectorcone2}.

Applying similar movements to each vertex of  $\mathrm{conv}(S)$ as many times as necessary we achieve the objective of Step 3.

\bigskip

Applying in a similar way Step 4, Step 5, etc., we hold the following after a finite number of steps:
\begin{itemize}
	\item The distances among two vertices of $\mathrm{conv}(S)$ are not involved in any equality of distances.
    \item $\|p_{i}-p_h\|\neq\|p_{j}-p_h\|$ $\forall i,j\in \{1,2,...,n\}, i\neq j$, and  $\forall h\in\{1,...,k\}$.
\end{itemize}

\smallskip
 Let us consider $p_h,p_r,p_l,u,q,p\in S$ that are ordered  in clockwise direction, such that $p_h, p_r, p_l$ are vertices of $\mathrm{conv}(S)$ and $u,q,p$ are not. After the previous steps, other equalities implicating vertices
of $\mathrm{conv}(S)$  could appear. Without loss of generality these equalities can be identified to one of the following:
\begin{itemize}
\item $\|p-p_r\|=\|p-p_l\|$
\item $\|p-p_r\|=\|q-p_l\|$
\item $\|p-p_r\|=\|q-p_h\|$
\item $\|p-p_r\|=\|q-u\|$.
\end{itemize}
If some of these equalities appears, reevaluate $m$ and move the vertex $p_r$ to
$$p_{r}+\frac{\delta}{2n}\left(\frac{p_{r}-p}{\|p_{r}-p\|}\right),$$
 where $\delta< \min\{m,\epsilon\}$. This movement breaks the coincidence,
and finally the vertices of $\mathrm{conv}(S)$ are not involved in any equality of distances.

In order to break the equalities of distances implicating other points of $S$, consider $S'=S\setminus\{\text{vertices of }\mathrm{conv}(S)\}$ and apply the above approach to $S'$ and $\mathrm{conv}(S')$, although evaluating $m$ as the minimum difference among all distances among points of $S$.

The process ends after a finite number of steps when there is not any equality among distances. Each point has been moved at most $n-1$ times and at a distance less than $\epsilon$ at the end.
\end{proof}

A positive answer to the following conjecture would allow to apply Algorithm \ref{maximum} directly to any set $S$ of $n$ points, even if the there is some equalities among the distances by considering the order $<_B$ for the distances.
\begin{figure}
	\begin{center}
		\scalebox{0.25}{\includegraphics{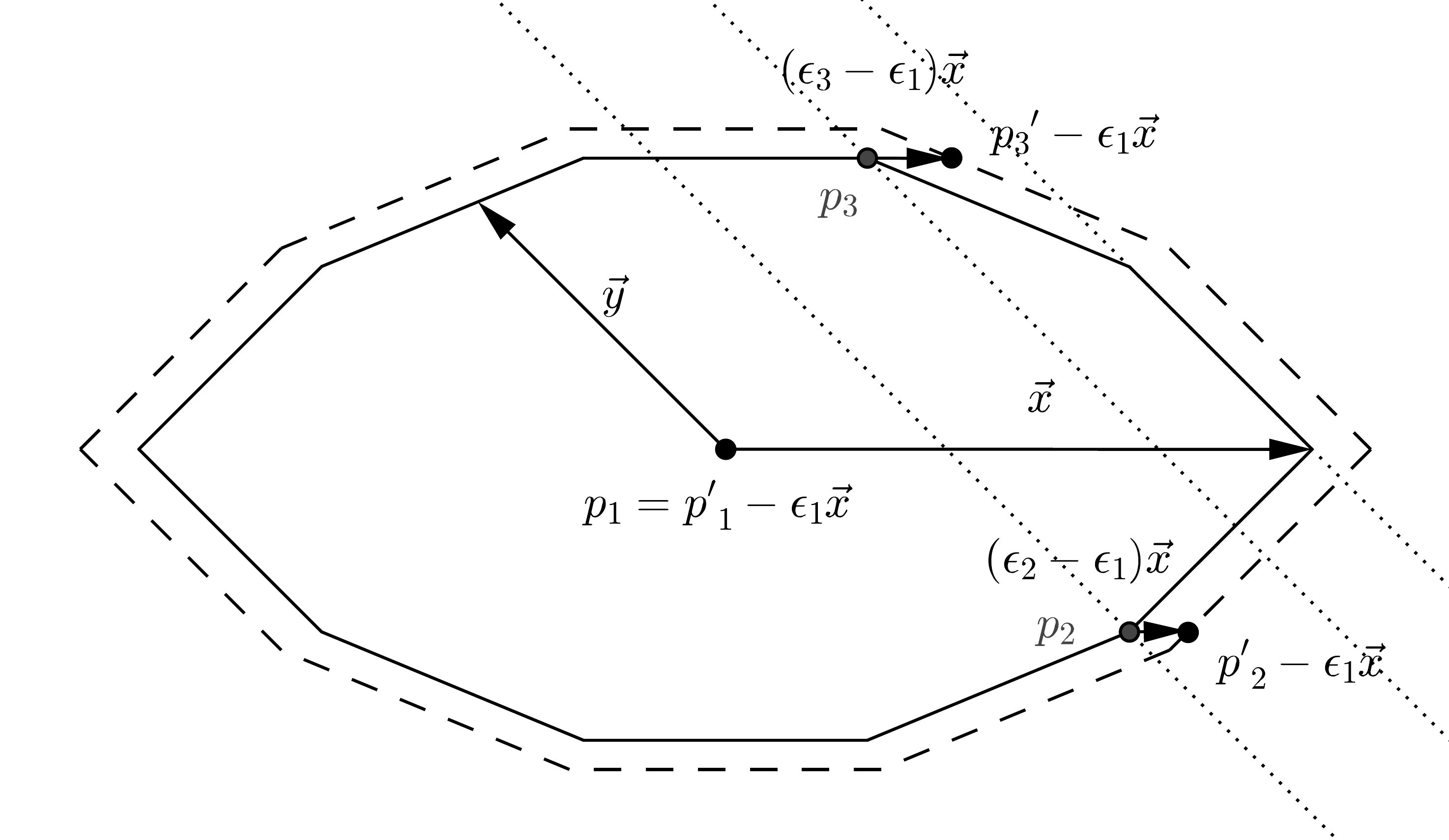}}\\
		\caption{\label{translation}$\|p_1-p_2\|=\|p_1-p_3\|=1$, $\|p'_1-p'_2\|=\|p'_1-p'_3\|=1.1$, with $p_i'=p_i+\epsilon_i$ ($i=1,2,3)$ and $\epsilon_1<\epsilon_2<\epsilon_3$.}.
	\end{center}
\end{figure}
\begin{conjecture}\label{conjetura}
	Let $\mathbb{M}^2$ be a strictly normed plane and $S$ be a finite set of $n$ points in $\mathbb{M}^2$. It is possible to denote $\{p_1,p_2,\cdots,p_n\}$ to the points of $S$ and order the distances among the points by some criterion $<_B$ such that
	         $\|p_i-p_j\| < \|p_k-p_l\|$ implies $\|p_i-p_j\| <_B \|p_k-p_l\|$, and
	for any $\epsilon>0$ there is a point set $S'=\{p'_1,p'_2,...,p'_n\}$ in a bijective correspondence with $S$. This correspondence verifies the following:
		\begin{enumerate}
        \item Each point in $S'$ is within distance $\epsilon$ of its corresponding point in $S$,
		\item The distances between points of $S'$ are all distinct,
		\item $\|p'_i-p'_j\|< \|p'_k-p'_l\|$ if and only if $\|p_i-p_j\| <_B \|p_k-p_l\|$.
	\end{enumerate}
\end{conjecture}


\bigskip

 Neither the statement of Conjecture 
 nor the one of Proposition \ref{movingpoints} are true for non strictly convex normed planes because there are bisectors that contain non empty interior, thus a slight movement of the points do not assure to break the coincidences. For instance, $p=(0,1)$ is an interior point of the bisector between $u=(-\frac{1}{2},0)$ and $v=(\frac{1}{2},0)$ with the maximum norm.

Let $x,y\in \mathbb{M}^2$ be a base such that $x$ is Birkhoff orthogonal to $y$ and $\|x\|=\|y\|=1$. The statement of  Conjecture 
is proved in \cite{MPSY} for $\mathbb{E}^2$ ordering the points lexicographically by their coordinates respect to $\{x,y\}$, and using the lexicographically order ($\|p_i-p_j\|, j, - i$) as criterion $<_B$ for the distances.
Under these conditions, for every $\epsilon>0$ there exist $0<\epsilon_1<\epsilon_2<\cdots<\epsilon_n<\epsilon$ such that  the movement from $p_i=(x_i,y_i)$ to $p_i'=(x_i+\epsilon_i,y_i)$ guarantees that:
\begin{enumerate}
	\item The point set $S'=\{p_1',p_2',\cdots,p_n'\}$ is in a bijective correspondence with $S$ (such that $p_i\to p_i'$).
	\item Each point in $S'$ is within distance $\epsilon$ of its corresponding point in $S$.
	\item The distances between points of $S'$ are all distinct.
	\item The order $<_B$ 
preserves the corresponding original order of distances.
\end{enumerate}

Nevertheless,
this approach does not work in every normed plane: in Figure \ref{translation}  the normed sphere is a polygon with ten sides
such that $\|p_1-p_2\|=\|p_1-p_3\|=1$ but the translation $p_i'=p_i+\epsilon_i$ ($i=1,2,3)$ does not break the equalities. Obviously, this example can be rounded in order to obtain a similar one but with a differentiable and strictly convex sphere. Our Proposition \ref{movingpoints} presents a partial positive answer to Conjecture 
for strictly convex normed planes.

\section*{Acknowledgments}
The authors gratefully thank Professor Pier Luigi Papini for his useful comments and suggestions.
\bibliographystyle{amsplain}

\begin{thebibliography}{00}
	
	
	
	
	\bibitem{Au} F. Aurenhammer, \textit{Voronoi diagrams - a survery of a fundamental geometric data structure}, ACM Computing Surveys (CSUR) \textbf{23} (1991),  345-405.
	
	
	\bibitem{Asano} T. Asano, B. Bhattacharya, M. Keil, and F. Yao, \textit{Clustering algorithms based on minimum and maximum spanning trees}, in Proc. 4th ACM Symposium on Computational Geometry, (1988),  252-257.
	
	

	
	\bibitem{E-M-P-R-W-W} H. Edelsbrunner, H.A. Mauer, F.P. Preparata, A.L. Rosenberg, E. Welzl, and D. Wood, \textit{Stabbing line segments},  BIT \textbf{22} (1982), 274-281.
	
	
	\bibitem{Gabow} H. Gabow, J. Bentley, and R. Tarjan, \textit{Scaling and related techniques for geometric problems}, in Proc. 16th ACM Symposium on Theory of Computing, (1984),  135-143.
	
	
	
	\bibitem{Hagauer-Rote} J. Hagauer, and G. Rote, \textit{Three-clustering of points in the plane}, Comput. Geom. \textbf{8} (1997), 87-95.
	
	


	\bibitem{JS} T. Jahn, \textit{On bisectors in normed planes}, Contributions to Discrete Mathematics, \textbf{10} (2015), 1-9.

\bibitem{Klein} R. Klein, \textit{Concrete and Abstract Voronoi Diagrams. Lecture Notes in Computer Science 400}, Springer-Verlag, New York, 1989.	


\bibitem{KMM} R. Klein, K. Mehlhorn, and S. Meiser,  \textit{Randomized Incremental Construction of Abstract Voronoi Diagrams}, Computational Geometry: Theory and Applications \textbf{3} (1993),  157-184.

	\bibitem{Lee} D.T. Lee, \textit{Farthest neighbor Voronoi diagrams and applications}, North-western University, Department of Electrical Engineering and Computer Science, 1980.
	

	\bibitem{Martin-Ya} P. Mart\'{\i}n, and D. Y\'{a}nez, \textit{Geometric clustering in normed planes}, Computational Geometry: Theory and Applications 78, \textbf{78} (2019),  50-60.

	
	
	\bibitem{MSW} H. Martini, K.J. Swanepoel, and G. Weiss, \textit{The geometry of Minkowski spaces -  a survey, Part I}, Expositiones Math. \textbf{19} (2001), 97-142.
	
	
		
	\bibitem{Mehlhorn} K. Mehlhorn, S. Meiser, and R. Rasch,\textit{ Furthest site abstract Vornoi diagrams}, International Journal of Computational Geometry and Applications, \textbf{11} (2001), 583-616.
	
	\bibitem{MPSY} C. Monma, M. Paterson, S. Suri, and F. Yao,\textit{  Computing Euclidean maximum spanning trees}, Algorithmica \textbf{5} (1990), 407-419.
	
	\bibitem{P-S} F.P. Preparata, and M.I. Shamos,  \emph{Computational Geometry}, Springer-Verlag, New York, 1985.
	
	\bibitem{Shamos-Hoey} M.I. Shamos and D. Hoey, \textit{Closest-point problems}, in Fundations of Computer Science, 1975, 16th Annual Symposium on, (1975),  151-162.
	
	
	
\end{thebibliography}

\end{document}